\documentclass[11pt]{article}
\usepackage{amsmath,epsfig,latexsym,amssymb,cite,graphicx,a4}
\newlength{\indentedwidth}
\newdimen\mathindent
\indentedwidth=\mathindent

\usepackage{epsf,graphics}
\DeclareMathAlphabet{\mathpzc}{OT1}{pzc}{m}{it}  
\begin{document}
\vskip 0.5cm
\begin{center}
{\Large \bf Eigenvalues of Casimir invariants for $U_q[osp(m|n)]$} 
\end{center}
\vskip 0.8cm
\centerline{K.A. Dancer%
\footnote{\tt dancer@maths.uq.edu.au}, 
M.D. Gould%
\footnote{\tt mdg@maths.uq.edu.au}
and J. Links%
\footnote{\tt jrl@maths.uq.edu.au}}
\vskip 0.9cm
\centerline{\sl\small 
 Centre for Mathematical Physics, School of Physical Sciences, }
\centerline{\sl\small The University of Queensland, Brisbane 4072,
 Australia.} 
\vskip 0.9cm
\begin{abstract}
\vskip0.15cm
\noindent
For each quantum superalgebra $U_q[osp(m|n)]$ with $m>2$, 
an infinite family
of Casimir invariants
is constructed. 
This is achieved by using an explicit form for the Lax operator. 
The eigenvalue of each Casimir invariant
on an arbitrary irreducible highest weight
module is also calculated. 
\end{abstract}

\setcounter{footnote}{0}
\def\thefootnote{\fnsymbol{footnote}}

\newenvironment{fig}{\linespread{1.0} \begin{figure}}{\end{figure}%
\linespread{1.3}}
\newcommand{\fl}{\hspace*{-\mathindent}}
\newcommand\phup{^{\phantom p}}
\newtheorem{definition}{Definition}[section]
\newtheorem{lemma}{Lemma}[section]
\newtheorem{theorem}{Theorem}[section]

\section{Introduction}
\noindent Representations of quantum superalgebras are known to
provide solutions to the Yang--Baxter equation and represent the
symmetries that underly supersymmetric exactly solvable (or integrable) models. Many
such examples have arisen in the context of modelling systems 
of strongly correlated electrons 
\cite{Foerster, Gonzalez, uq, uqagain, Martins}.
More recently, the properties of solvability and supersymmetry have been applied
to other areas, such as the solution of the Kondo 
model \cite{bk}, integrable superconformal field theory \cite{kz} and disordered systems \cite{efs}. 
Developing the representation theory of the quantum
superalgebras is a useful step towards the complete understanding of such models.
However, in many respects the representation theory of quantum superalgebras is not a 
straightforward generalisation of the quantum algebra case, principally because not all 
representations of quantum superalgebras are unitary \cite{gs}. This leads, for example, to
the existence of indecomposable representations not arising in   
the quantum algebra case, which generally make the analysis of supersymmetric models
problematic (e.g. see \cite{efs}).

In this paper we construct the Casimir invariants (central elements) 
of quantised
orthosymplectic superalgebras.  
Our method of construction 
follows from the general results of \cite{ZhangGould}
and the
explicit form of the Lax operator obtained in \cite{Lax}.
A fundamental 
problem  
is to determine the eigenvalues of the Casimir invariants when acting 
on an arbitrary finite-dimensional irreducible module. To date, 
the eigenvalues have only been calculated for 
the type I quantum superalgebras \cite{GLZ, LinksZhang}, while the results for 
$U_q[osp(1|n)]$ follow from an isomorphism derived in \cite{Zhang}.  
In this paper we perform the
calculations for the remaining non-exceptional quantum superalgebras,
namely $U_q[osp(m|n)]$ for $m >2$. 
The procedure we use for calculating the eigenvalues of the Casimir invariants 
when acting on 
any irreducible module is based on the early work by Perelomov and Popov 
\cite{Perelomov, Popov} and Nwachuku and Rashid \cite{Nwachuku}.  In
doing so we follow the method used in \cite{Bincer} and
\cite{Scheunert83} for the classical general and orthosymplectic
superalgebras respectively, which was adapted in \cite{LinksZhang} to
cover $U_q[gl(m|n)]$.  Although the concepts are much the same as in
those cases, the combination of the $q$-deformation and the more
complex root system of $U_q[osp(m|n)]$ makes the calculations in this
paper more technically challenging.

In the following section we introduce our notation for $U_q[osp(m|n)]$
and state the Lax operator.  In Section 3 we develop the formulae for the
Casimir invariants of $U_q[osp(m|n)]$.  The bulk of the calculations
are in Section 4, where the eigenvalues of the Casimir invariants are
derived in detail.
\section{The quantised orthosymplectic superalgebra $U_q[osp(m|n)]$}

\noindent The quantum superalgebra $U_q[osp(m|n)]$ is a
$q$-deformation of the classical orthosymplectic superalgebra.  A
brief explanation of $U_q[osp(m|n)]$ is given below, with a more
thorough introduction to $osp(m|n)$ and the $q$-deformation to be
found in \cite{Lax}.

First we need to define the notation.  The grading of $a$ is
denoted by $[a]$, where

\begin{equation*}
[a] =
   \begin{cases}
     0, \qquad \; a=i, & 1 \leq i \leq m, \\ 1, \qquad \; a = \mu, & 1
     \leq \mu \leq n.
   \end {cases}
\end{equation*}

\noindent Throughout this paper we use Greek letters $\mu, \nu$
etc. to denote odd indices and Latin letters $i, j$ etc. for even
indices.  If the grading is unknown, the usual $a, b, c$ etc. are
used.  Which convention applies will be clear from the context. 
 Throughout the paper we also use the symbols $\overline{a}$
and $\xi_a$, which are given by:

\begin{equation*}
\overline{a}=
  \begin{cases}
     m + 1 - a, & [a]=0, \\ n + 1 - a, & [a]=1,
  \end {cases}
\qquad \text{and} \quad \xi_{a} =
  \begin{cases}
     1, & [a] = 0, \\ (-1)^a, & [a]=1.
  \end{cases}
\end{equation*}

As a weight system for $U_q[osp(m|n)]$ we take the set $\{
\varepsilon_i,\; 1 \leq i \leq m \} \cup \{ \delta_\mu,\; 1 \leq \mu
\leq n \}$, where \mbox{$\varepsilon_{\overline{i}} = -
\varepsilon_i$} and \mbox{$\delta_{\overline{\mu}} = - \delta_\mu$}.
Conveniently, when $m=2l+1$ this implies $\varepsilon_{l+1} =
-\varepsilon_{l+1} = 0$.  Acting on these weights, we have the
invariant bilinear form defined by:

\begin{equation*}
(\varepsilon_i, \varepsilon_j) = \delta^i_j, \quad (\delta_\mu,
\delta_\nu) = -\delta^\mu_\nu, \quad (\varepsilon_i, \delta_\mu) = 0,
\qquad 1 \leq i, j \leq l, \quad 1 \leq \mu, \nu \leq k.
\end{equation*}

\noindent When describing an object with unknown grading indexed by
$a$ the weight will be described generically as $\varepsilon_a$.  This
should not be assumed to be an even weight.

The even positive roots of $U_q[osp(m|n)]$ are composed
entirely of the usual positive roots of $o(m)$ together with those of
$sp(n)$, namely:

\begin{alignat*}{3}
&\varepsilon_i \pm \varepsilon_j, & \qquad & 1 \leq i < j \leq l, \\
&\varepsilon_i, && 1 \leq i \leq l &&\quad \text{when } m=2l+1, \\
&\delta_\mu + \delta_\nu, && 1 \leq \mu,\nu \leq k, \\ &\delta_\mu -
\delta_\nu, && 1 \leq \mu < \nu \leq k.
\end{alignat*}

\noindent The root system also contains a set of odd positive roots,
which are:

\begin{equation*}
\delta_\mu + \varepsilon_i, \qquad 1 \leq \mu \leq k,\;1 \leq i \leq
m.
\hspace{2cm}
\end{equation*}

\noindent Throughout this paper we choose to use the following set of
simple roots:

\begin{align*}
&\alpha_i = \varepsilon_i - \varepsilon_{i+1}, \hspace{11mm} 1 \leq i
  < l, \notag \\ &\alpha_l =
\begin{cases} 
\varepsilon_{l} + \varepsilon_{l-1},\quad & m=2l, \\ \varepsilon_l,
&m=2l+1,
\end{cases} \notag \\
&\alpha_\mu = \delta_\mu - \delta_{\mu+1}, \hspace{9mm} 1 \leq \mu <
k,\notag\\ &\alpha_s = \delta_k - \varepsilon_1.
\end{align*}

\noindent Note this choice is only valid for $m >2$.

In $U_q[osp(m|n)]$ the graded commutator is realised by

\begin{equation*}
[A,B] = AB - (-1)^{[A][B]} BA
\end{equation*}
 
\noindent and tensor product multiplication is given by

\begin{equation*}
(A \otimes B) (C \otimes D) = (-1)^{[B][C]} (AC \otimes BD).
\end{equation*}

\noindent Using these conventions, we have:

\begin{definition} \label{def} The quantum superalgebra $U_q[osp(m|n)]$ is generated by simple 
generators $e_a, f_a, h_a$ subject to the relations:

\begin{alignat}{2}
&[h_a, e_b] = (\alpha_a, \alpha_b) e_b, && \notag \\ &[h_a, f_b] = -
(\alpha_a, \alpha_b) f_b, && \notag \\ &[h_a, h_b] = 0,&& \notag \\
&[e_a, f_b] = \delta^a_b \frac{(q^{h_a} - q^{-h_a})}{(q - q^{-1})},&&
\notag \\ &[e_a, e_a] = [f_a,f_a]=0 & \quad &\text{for }\; (\alpha_a,
\alpha_a)=0,\notag
\end{alignat}
We remark that $U_q[osp(m|n)]$ has the structure of a quasi-triangular
Hopf superalgebra. In 
particular, there is a linear mapping known as the \textit{coproduct},
$\Delta: U_q[osp(m|n)] \rightarrow U_q[osp(m|n)]^{\otimes
2}$, which is defined on the simple generators by:

\begin{align*}
&\Delta (e_a) = q^{\frac{1}{2}h_a} \otimes e_a + e_a \otimes
  q^{-\frac{1}{2} h_a}, \notag \\ &\Delta (f_a) = q^{\frac{1}{2}h_a}
  \otimes f_a + f_a \otimes q^{-\frac{1}{2} h_a}, \\ 
&\Delta (q^{\pm \frac{1}{2}h_a}) = q^{\pm \frac{1}{2}h_a} \otimes q^{\pm
  \frac{1}{2}h_a}, 
\end{align*}

\noindent and extends to arbitrary elements according to the homomorphism 
 property, namely:
\begin{align*}
\Delta (AB) = \Delta(A) \Delta(B). \label{coprod2}
\end{align*}

\end{definition}

\noindent There are further defining relations such as the $q$-Serre
relations, but they are not needed in this paper.

The quasi-triangular property guarantees the existence of
a universal $R$-matrix, which provides a solution to the 
Yang--Baxter equation.  Before elaborating, we need to introduce the
graded twist map.

 The \textit{graded twist map} $T:U_q[osp(m|n)]^{\otimes 2}
\rightarrow U_q[osp(m|n)]^{\otimes 2}$ is given by

\begin{equation*}
T(a \otimes b) = (-1)^{[a][b]} (b \otimes a).
\end{equation*}

\noindent For convenience $T \circ \Delta$, the twist map composed
with the coproduct, is denoted $\Delta^T$.  Then a \textit{universal
$R$-matrix}, $\mathcal{R}$, is an even, non-singular element of
$U_q[osp(m|n)]^{\otimes 2}$ satisfying the following properties:

\begin{align}
&\mathcal{R} \Delta (a) = \Delta^T (a)\mathcal{R}, \quad \forall a \in
  U_q[osp(m|n)], \notag \\ &(\text{id} \otimes \Delta) \mathcal{R} =
  \mathcal{R}_{13} \mathcal{R}_{12}, \notag \\ &(\Delta \otimes
  \text{id}) \mathcal{R} = \mathcal{R}_{13} \mathcal{R}_{23}.
  \label{Requations}
\end{align}

\noindent Here $\mathcal{R}_{ab}$ represents a copy of $\mathcal{R}$
acting on the $a$ and $b$ components respectively of $U_1 \otimes U_2
\otimes U_3$, where each $U$ is a copy of the quantum superalgebra
$U_q[osp(m|n)]$.  When $a>b$ the usual grading term from the twist map
is included, so for example $\mathcal{R}_{21} = [\mathcal{R}^T]_{12}
$, where $\mathcal{R}^T = T (\mathcal{R})$ is the \textit{opposite
universal $R$-matrix}.

The $R$-matrix is significant because it is a solution to the 
Yang--Baxter equation, which is prominent in the study of integrable
systems \cite{Baxter}:

\begin{equation*}
\mathcal{R}_{12} \mathcal{R}_{13} \mathcal{R}_{23} = \mathcal{R}_{23}
  \mathcal{R}_{13} \mathcal{R}_{12}
\end{equation*}

\noindent A superalgebra may contain many different universal
$R$-matrices, but there is always a unique one belonging to
$U_q[osp(m|n)]^- \otimes U_q[osp(m|n)]^+$, with its opposite $R$-matrix
in $U_q[osp(m|n)]^+ \otimes U_q[osp(m|n)]^-$.  Here $U_q[osp(m|n)]^-$
is the Hopf subsuperalgebra generated by the lowering generators $\{f_a\}$ and
Cartan elements $\{h_a\}$, while $U_q[osp(m|n)]^+$ is generated by the raising
generators $\{e_a\}$ and the Cartan elements.  These particular $R$-matrices arise
out of the $\mathbb{Z}_2$-graded version of Drinfeld's double construction 
\cite{gzb}.  
In this paper
we consider the universal $R$-matrix belonging to $U_q[osp(m|n)]^-
\otimes U_q[osp(m|n)]^+$.

 We also need to define the vector representation for
$U_q[osp(m|n)]$. Let $\text{End} \; V$ be the space of endomorphisms of
$V$, an $(m+n)$-dimensional vector space.  Then the irreducible
\textit{vector representation} $\pi: U_q[osp(m|n)] \rightarrow
\text{End} \; V$ acts on the $U_q[osp(m|n)]$ generators as 
given in Table \ref{vector}, where $E^a_b$ is the elementary matrix with a 1 
in the $(a,b)$ position and zeroes elsewhere.

One quantity that repeatedly arises in calculations for both
classical and quantum Lie superalgebras is $\rho$, the \textit{graded
half-sum of positive roots}.  In the case of $U_q[osp(m|n)]$ it is
given by:

\begin{equation*}
\rho = \frac{1}{2} \sum_{i=1}^l (m-2i) \varepsilon_i + \frac{1}{2}
\sum_{\mu=1}^k (n-m+2-2\mu) \delta_\mu.
\end{equation*}

\noindent This satisfies the property $(\rho, \alpha) = \frac{1}{2}
(\alpha, \alpha)$ for all simple roots $\alpha$.

\begin{table}[ht]
\caption{The action of the vector representation $\pi$ on the simple 
generators of $U_q[osp(m|n)]$} 
\label{vector}
\centering
\begin{tabular}{|l|l|l|l|}\hline
\multicolumn{1}{|c|}{$\alpha_a$}& \multicolumn{1}{c|}{$\pi(e_a)$} & 
  \multicolumn{1}{c|}{$\pi(f_a)$} & \multicolumn{1}{c|}{$\pi(h_a)$} \\ \hline 
$\alpha_i, 1 \leq i < l$ &
  $E^i_{i+1} - E^{\overline{i+1}}_{\overline{i}}$&
  $E_i^{i+1} - E_{\overline{i+1}}^{\overline{i}}$&
  $E^i_i - E^{\overline{i}}_{\overline{i}} - E^{i+1}_{i+1}
             + E^{\overline{i+1}}_{\overline{i+1}}$\\ 
$\alpha_l, \, m=2l$ &
  $E^{l-1}_{\overline{l}} - E^l_{\overline{l-1}}$ &
  $E_{l-1}^{\overline{l}} - E_l^{\overline{l-1}}$&
  $E^{l-1}_{l-1} + E^l_l - E^{\overline{l-1}}_{\overline{l-1}}
             - E^{\overline{l}}_{\overline{l}}$ \\
$\alpha_l, \, m=2l+1$ 
  &$E^l_{l+1} - E^{l+1}_{\overline{l}} $ &
  $E_l^{l+1} - E_{l+1}^{\overline{l}}$&
  $E^l_l - E^{\overline{l}}_{\overline{l}} $\\
$\alpha_\mu, 1 \leq \mu <k$ & 
  $E^\mu_{\mu+1} + E^{\overline{\mu+1}}_{\overline{\mu}}$&
  $E_\mu^{\mu+1} + E_{\overline{\mu+1}}^{\overline{\mu}}$&
  $E^{\mu+1}_{\mu+1} - E^{\overline{\mu+1}}_{\overline{\mu+1}}
               - E^\mu_\mu + E^{\overline{\mu}}_{\overline{\mu}}$ \\
$\alpha_s$ &
  $E^{\mu=k}_{i=1} + (-1)^k E^{\overline{i=1}}_{\overline{\mu=k}}$&
  $- E^{i=1}_{\mu=k} + (-1)^k E_{\overline{i=1}}^{\overline{\mu=k}}$&
  $- E^{i=1}_{i=1} + E^{\overline{i}=\overline{1}}_{\overline{i}=\overline{1}} 
   -  E^{\mu=k}_{\mu=k} + E^{\overline{\mu}=\overline{k}}_
           {\overline{\mu}=\overline{k}}$ \\ \hline
\end{tabular}
\end{table}

\subsection{The Lax operator for $U_q[osp(m|n)]$}

\noindent Let $\mathcal{R}$ be the universal $R$-matrix of
$U_{q}[osp(m|n)]$ and $\pi$ the vector representation.  The Lax
operator associated with $\mathcal{R}$ is given by

\begin{equation*}
R = (\pi \otimes \text{id}) \mathcal{R} \in (\text{End} \; V) \otimes
U_{q}[osp(m|n)].
\end{equation*}

\noindent It has been shown in \cite{Lax} that the Lax operator is
given by

\begin{equation*}
R = \sum_a E^a_a \otimes q^{h_{\varepsilon_a}} + (q-q^{-1})
  \sum_{\varepsilon_a < \varepsilon_b} (-1)^{[b]} E^a_b \otimes
  q^{h_{\varepsilon_a}} \hat{\sigma}_{ba},
\end{equation*}

\noindent where the simple operators $\hat{\sigma}_{ba}$ are given by:

\begin{alignat*}{2}
&\hat{\sigma}_{i\, i+1} = - \hat{\sigma}_{\overline{i+1}\,
  \overline{i}} = q^{\frac{1}{2}} e_i q^{\frac{1}{2} h_i},&&1 \leq i <
  l, \notag \\ &\hat{\sigma}_{l-1\, \overline{l}} = -
  \hat{\sigma}_{l\, \overline{l-1}} = q^{\frac{1}{2}} e_l
  q^{\frac{1}{2} h_l}, && m=2l, \notag \\ &\hat{\sigma}_{l\,
  \overline{l}} = 0, && m=2l, \notag \\ &\hat{\sigma}_{l\, l+1} = -
  q^{-\frac{1}{2}} \hat{\sigma}_{l+1\, \overline{l}} = e_l
  q^{\frac{1}{2} h_l}, && m=2l+1, \notag \\ &\hat{\sigma}_{\mu \, \mu
  +1} = \hat{\sigma}_{\overline{\mu+1}\,\overline{\mu}} =
  q^{-\frac{1}{2}} e_\mu q^{\frac{1}{2} h_\mu}, &\qquad&1 \leq \mu <k,
  \notag \\ &\hat{\sigma}_{\mu=k \, i=1} = (-1)^k q\, \hat{\sigma}_{i
  = \overline{1} \, \overline{\mu} = \overline{k}} = q^{\frac{1}{2}}
  e_s q^{\frac{1}{2} h_s};
\end{alignat*}

\noindent and the remaining operators can be calculated using

\noindent (i) the $q$-commutation relations

\begin{equation*}
q^{(\alpha_c, \varepsilon_b)} \hat{\sigma}_{ba} e_c q^{\frac{1}{2}
  h_c} - (-1)^{([a]+[b])[c]} q^{-(\alpha_c, \varepsilon_a)} e_c
  q^{\frac{1}{2}h_c} \hat{\sigma}_{ba} = 0, \quad \varepsilon_b >
  \varepsilon_a,
\end{equation*}
      
\noindent where neither $\varepsilon_a - \alpha_c$ nor $\varepsilon_b
+ \alpha_c$ equals any $\varepsilon_x$; and

\noindent (ii) the induction relations

\begin{equation*}
\hat{\sigma}_{ba} = q^{-(\varepsilon_b,\varepsilon_a)}
  \hat{\sigma}_{bc} \hat{\sigma}_{ca} - q^{-(\varepsilon_c,
  \varepsilon_c)} (-1)^{([b]+[c]) ([a]+[c])} \hat{\sigma}_{ca}
  \hat{\sigma}_{bc}, \quad \varepsilon_b > \varepsilon_c >
  \varepsilon_a,
\end{equation*}

\noindent where $ c \neq \overline{b} \text{ or } \overline{a}$.

To define the opposite Lax operator $R^T = (\pi \otimes
\text{id}) \mathcal{R}$ we require the graded conjugation action
$\dagger$, which is defined on the simple generators by (see \cite{Lax}):

\begin{equation*}
e_a^\dagger = f_a, \qquad f_a^\dagger = (-1)^{[a]} e_a, \qquad
h_a^\dagger = h_a.
\end{equation*}

\noindent It is consistent with the coproduct and extends naturally to
all remaining elements of $U_q[osp(m|n)]$, satisfying the properties:

\begin{align*}
&(\sigma^a_b)^\dagger = (-1)^{[a]([a]+[b])} \sigma^b_a,\\
&(ab)^\dagger = (-1)^{[a][b]} b^\dagger a^\dagger, \\ &(a \otimes
b)^\dagger = a^\dagger \otimes b^\dagger, \\ &\Delta (a)^\dagger =
\Delta (a^\dagger).
\end{align*}

\noindent Then the opposite $R$-matrix is given by

\begin{equation*}
R^T = \sum_a E^a_a \otimes q^{h_{\varepsilon_a}} + (q-q^{-1})
\sum_{\varepsilon_b > \varepsilon_a} (-1)^{[a]} E^b_a \otimes
\hat{\sigma}_{ab} q^{h_{\varepsilon_a}},
\end{equation*}

\noindent where

\begin{equation*}
\hat{\sigma}_{ab} = (-1)^{[b]([a]+[b])} \hat{\sigma}_{ba}^\dagger,
\quad \varepsilon_b > \varepsilon_a.
\end{equation*}

\section{Casimir Invariants of $U_q[osp(m|n)]$}

\noindent We now use the Lax operator to construct a family of Casimir
invariants and then to calculate their eigenvalues when acting on an
irreducible highest weight module.  
Before constructing the Casimir invariants, however, we need to define
a new object. Let $h_\rho$ be the unique element of the Cartan
subalgebra H satisfying

\begin{equation*}
\alpha_i(h_\rho) = (\rho, \alpha_i), \qquad \forall \alpha_i \in H^*.
\end{equation*}

\noindent Then from \cite{ZhangGould} we have the following theorem:

\begin{theorem}
Let $\mathcal{V}$ be the representation space of $\tau$, an arbitrary 
finite-dimensional representation of $U_q[osp(m|n)]$.  If $\Gamma \in
(\text{End } \mathcal{V}) \otimes U_q[osp(m|n)]$ satisfies

\begin{equation} \label{partial}
\partial (a) \Gamma = \Gamma \partial (a), \quad \forall a \in
U_q[osp(m|n)],
\end{equation}

\noindent where $\partial \equiv (\pi \otimes \text{id}) \Delta$, then

\begin{equation*}
C = (str \otimes \mathrm{id}) (\tau (q^{2h_\rho}) \otimes I) \Gamma
\end{equation*}

\noindent belongs to the centre of $U_q[osp(m|n)]$.  Above $str$
denotes the supertrace.
\end{theorem}

Now choose $\tau$ to be the vector representation $\pi$.  Recalling
that the universal $R$-matrix satisfies

\begin{equation*} 
\mathcal{R} \Delta(a) = \Delta^T(a) \mathcal{R}, \qquad \forall a \in
U_q[osp(m|n)],
\end{equation*}

\noindent it is clear that

\begin{equation*}
\partial(a) R^T R = R^T R \,\partial(a), \qquad \forall a \in
U_q[osp(m|n)].
\end{equation*}

\noindent Hence if we set $A \in (\text{End}\;V) \otimes
U_q[osp(m|n)]$ to be

\begin{equation*}
A = \frac{(R^T R - I \otimes I)}{(q-q^{-1})},\
\end{equation*}

\noindent the operators $A^l$ will satisfy condition \eqref{partial}
for all non-negative integers $l$. Thus the operators $C_l$ defined as

\begin{equation*}
C_l = (str \otimes \text{id})(\pi(q^{2h_p}) \otimes I) A^l, \qquad l
\in \mathbb{Z}^+,
\end{equation*}

\noindent form a family of Casimir invariants.  Here $A$ coincides
with the matrix of Jarvis and Green \cite{Jarvis} in the classical
limit $q \rightarrow 1$, as do the invariants $C_l$.
 
Now write the Lax operator $R$ and its opposite $R^T$ in the
form

\begin{align*}
R &= I \otimes I + (q-q^{-1}) \sum_{\varepsilon_b \geq \varepsilon_a}
  E^a_b \otimes X^b_a, \\ R^T &= I \otimes I + (q-q^{-1})
  \sum_{\varepsilon_b \leq \varepsilon_a} E^a_b \otimes X^b_a.
\end{align*}

\noindent In terms of the operators $\hat{\sigma}_{ba}$, this implies

\begin{equation*}
X^b_a = \begin{cases} \frac{q^{h_{\varepsilon_a}}-I}{q-q^{-1}}, &a=b,
\\ (-1)^{[b]} q^{h_{\varepsilon_a}} \hat{\sigma}_{ba}, \quad
&\varepsilon_a < \varepsilon_b, \\ (-1)^{[b]} \hat{\sigma}_{ba}
q^{h_{\varepsilon_b}}, & \varepsilon_a > \varepsilon_b.
\end{cases}
\end{equation*}

\noindent Writing $A$ as

\begin{equation*}
A = \sum_{a,b} E^a_b \otimes A^b_a,
\end{equation*}

\noindent we obtain

\begin{equation*}
A^b_a = (1+\delta^a_b) X^b_a + (q-q^{-1}) \sum_{\varepsilon_c \leq
\varepsilon_a, \varepsilon_b} (-1)^{([a]+[c])([b]+[c])} X^c_a X^b_c.
\end{equation*}

\noindent This produces a family of Casimir invariants

\begin{equation*}
C_l = \sum_a (-1)^{[a]} q^{(2\rho, \varepsilon_a)} {A^{(l)}}^a_a,
\end{equation*}

\noindent where the operators ${A^{(l)}}^b_a$ are recursively defined
as

\begin{equation} \label{Alba}
{A^{(l)}}^b_a = \sum_c (-1)^{([a]+[c])([b]+[c])} {A^{(l-1)}}^c_a
A^b_c.
\end{equation}

Note that $A$ corresponds to the matrix $A$ given for the non-graded case in
\cite{LinksGould91}.  Following a line of reasoning similar to that in
\cite{Gould87}, it can be shown that when acting on an irreducible module
$V(\Lambda)$, $A$ satisfies the following polynomial identity:

\begin{equation*}
\prod_{a=1}^{m+n} (A - \alpha_a(\Lambda)I) = 0
\end{equation*}

\noindent where

\begin{equation*}
\alpha_a(\Lambda) = \frac{q^{(\varepsilon_a, \varepsilon_a + 2\Lambda + 2\rho) - C(\Lambda_0)}-1}{q-q^{-1}}
\end{equation*}

\noindent and $C(\Lambda_0) = (\delta_1, \delta_1 + 2\rho) = m-n-1.$  In the limit $q \rightarrow 1$ this reduces to the identity given in \cite{Gould87}.


\section{Eigenvalues of the Casimir invariants}

\noindent Now that we have found a family of Casimir invariants, we
wish to calculate their eigenvalues on a general irreducible
finite-dimensional module. Let $V(\Lambda)$ be an arbitrary
irreducible finite-dimensional module with highest weight $\Lambda$
and highest weight state $| \Lambda \rangle$.  Define $t_a^{(l)}$ to
be the eigenvalue of ${A^{(l)}}^a_a$ on this state, so

\begin{equation*}
 {A^{(l)}}^a_a |\Lambda \rangle = t_a^{(l)} |\Lambda \rangle.
\end{equation*}

\noindent Once we have calculated $t_a^{(l)}$ we will use this result
to find the eigenvalues of the Casimir invariants $C_l$.

To evaluate $t_a^{(l)}$, note that if $\varepsilon_b >
\varepsilon_a$ then ${A^{(l)}}^b_a$ is a raising operator, implying
${A^{(l)}}^b_a |\Lambda \rangle = 0$. Thus from equation \eqref{Alba}
we deduce

\begin{align*}
t_a^{(l)} |\Lambda \rangle &= t_a^{(l-1)} t_a^{(1)} |\Lambda \rangle +
  \sum_{\varepsilon_a < \varepsilon_b} (-1)^{[a]+[b]} {A^{(l-1)}}^b_a
  A^a_b |\Lambda \rangle \notag \\ &= t_a^{(l-1)} t_a^{(1)} |\Lambda
  \rangle + \sum_{\varepsilon_a< \varepsilon_b} (-1)^{[a]+[b]}
  {A^{(l-1)}}^b_a \bigl[ X^a_b + (q-q^{-1}) X^a_b X^a_a \bigr]
  |\Lambda \rangle \notag \\ &= t_a^{(l-1)} t_a^{(1)} |\Lambda \rangle
  + \sum_{\varepsilon_a< \varepsilon_b} (-1)^{[a]+[b]} q^{(\Lambda,
  \varepsilon_a)} {A^{(l-1)}}^b_a X^a_b |\Lambda \rangle.
\end{align*}

\noindent Now we know that

\begin{equation} \label{AdelX}
A^l \partial (X^a_b) = \partial (X^a_b) A^l.
\end{equation}

\noindent This can be used to calculate ${A^{(l)}}^b_a X^a_b|\Lambda
\rangle$ for $\varepsilon_a< \varepsilon_b$. First we need an
expression for $\Delta(X^a_b)$. The $R$-matrix properties give

\begin{align*}
&(\Delta \otimes I) R = R_{13} R_{23} \\ \Rightarrow \quad &(I \otimes
\Delta) R^T = R^T_{12} R^T_{13}.
\end{align*}

\noindent In terms of $X^a_b$, this implies

\begin{align*}
I \otimes I \otimes I + &(q-q^{-1}) \sum_{\varepsilon_a \leq
  \varepsilon_b} E^b_a \otimes \Delta (X^a_b) \\ &= \bigl( I \otimes I
  \otimes I + (q-q^{-1}) \sum_{\varepsilon_a \leq \varepsilon_b} E^b_a
  \otimes X^a_b \otimes I \bigr) \\ & \hspace{2cm}\times \bigl( I
  \otimes I \otimes I + (q-q^{-1}) \sum_ {\varepsilon_a \leq
  \varepsilon_b} E^b_a \otimes I \otimes X^a_b \bigr) \\ &= I \otimes
  I \otimes I + (q-q^{-1}) \sum_{\varepsilon_a \leq \varepsilon_b}
  E^b_a \otimes (X^a_b \otimes I + I \otimes X^a_b) \\ & \hspace{2cm}
  + (q-q^{-1})^2 \sum_{\varepsilon_a \leq \varepsilon_c \leq
  \varepsilon_b} (-1)^{([a]+[c])([b]+[c])} E^b_a \otimes X^c_b \otimes
  X^a_c.
\end{align*}

\noindent Hence for all $\varepsilon_a < \varepsilon_b$

\begin{equation*}
 \Delta (X^a_b) = X^a_b \otimes I + I \otimes X^a_b + (q-q^{-1})
 \sum_{\varepsilon_a \leq \varepsilon_c \leq \varepsilon_b}
 (-1)^{([a]+[c])([b]+[c])} X^c_b \otimes X^a_c.
\end{equation*}

We also need an expression for $\pi (X^a_b)$ for
$\varepsilon_a \leq \varepsilon_b$.  In \cite{Lax} we found the
generators for $R^T$ in the vector representation are given by

\begin{equation*}
\hat{\sigma}_{ab} q^{h_{\varepsilon_a}} = E^a_b - (-1)^{[a]([a]+[b])}
\xi_a \xi_b q^{(\rho, \varepsilon_a - \varepsilon_b)}
E^{\overline{b}}_{\overline{a}}, \qquad \varepsilon_a < \varepsilon_b.
\end{equation*}

\noindent From this we deduce that

\begin{equation*}
\pi (X^a_b) = (-1)^{[a]} E^a_b - (-1)^{[a][b]} \xi_a \xi_b q^{(\rho,
  \varepsilon_a - \varepsilon_b)} E^{\overline{b}}_{\overline{a}},
  \qquad \varepsilon_a < \varepsilon_b.
\end{equation*}

\noindent Also, we know

\begin{align*}
\pi(X^a_a) &= (q-q^{-1})^{-1} \pi (q^{h_{\varepsilon_a}} - I) \\ &=
(q-q^{-1})^{-1} (q^{(\varepsilon_a, \varepsilon_a)(E^a_a -
E^{\overline{a}}_{\overline{a}})} -I).
\end{align*}

\noindent Applying these, we find that if $\varepsilon_a <
\varepsilon_b$ then

\begin{align*}
\partial (X^a_b) &= (\pi \otimes I) \Delta (X^a_b) \notag \\ &=
\pi(X^a_b) \otimes \bigl( I + (q-q^{-1}) X^a_a \bigr) + \bigl( I+
(q-q^{-1}) \pi (X^b_b) \bigr) \otimes X^a_b \notag \\ & \quad +
(q-q^{-1}) \sum_{\varepsilon_a < \varepsilon_c < \varepsilon_b}
(-1)^{([a]+[c])([b]+[c])} \pi(X^c_b) \otimes X^a_c \notag \\ &= \bigl(
(-1)^{[a]} E^a_b - (-1)^{[a][b]} \xi_a \xi_b q^{(\rho, \varepsilon_a -
\varepsilon_b)} E^{\overline{b}}_{\overline{a}} \bigr) \otimes
q^{h_{\varepsilon_a}} + q^{(\varepsilon_b, \varepsilon_b)(E^b_b -
E^{\overline{b}}_{\overline{b}})} \otimes X^a_b \notag \\ & \quad +
(q-q^{-1}) \sum_{\varepsilon_a < \varepsilon_c < \varepsilon_b}
(-1)^{([a]+[c])([b]+[c])} \notag \\ & \hspace{45mm}\times \bigl(
(-1)^{[c]} E^c_b - (-1)^{[b][c]} \xi_b \xi_c q^{(\rho, \varepsilon_c -
\varepsilon_b)} E^{\overline{b}}_{\overline{c}} \bigr) \otimes X^a_c.
\end{align*}

\noindent Substituting this expression into equation \eqref{AdelX} and
equating the $(a,b)$ entries, we find

\begin{align*}
(-1)^{[a]} &{A^{(l)}}^a_a q^{h_{\varepsilon_a}} -
  \delta^a_{\overline{b}} (-1)^{[a][b]} \xi_a \xi_b q^{(\rho,
  \varepsilon_a - \varepsilon_b)} {A^{(l)}} ^a_a q^{h_{\varepsilon_a}}
  + q^{(\varepsilon_b, \varepsilon_b)} {A^{(l)}}^b_a X^a_b \notag \\ &
  \qquad +(q-q^{-1}) \sum_{\varepsilon_a< \varepsilon_c
  <\varepsilon_b} \bigl( (-1)^{[c]} {A^{(l)}}^c_a X^a_c -
  \delta^b_{\overline{c}} (-1)^{[b][c]} \xi_b \xi_c q^{(\rho,
  \varepsilon_c - \varepsilon_b)} {A^{(l)}}^{\overline{b}}_a X^a_c
  \bigr) \notag \\ =& (-1)^{[a]} q^{h_{\varepsilon_a}} {A^{(l)}}^b_b -
  \delta^a_{\overline{b}} (-1)^{[a][b]} \xi_a \xi_b q^{(\rho,
  \varepsilon_a - \varepsilon_b)} q^{h_{\varepsilon_a}} {A^{(l)}}^b_b
  + (-1)^{[a]+[b]} q^{(\varepsilon_a, \varepsilon_b)} X^a_b
  {A^{(l)}}^b_a \notag \\ &\qquad -(q-q)^{-1} \delta^a_{\overline{b}}
  \sum_{\varepsilon_a<\varepsilon_c < \varepsilon_b} (-1)^{[b][c]}
  \xi_b \xi_c q^{(\rho, \varepsilon_c - \varepsilon_b)} X^a_c
  {A^{(l)}}^b_{\overline{c}}.
\end{align*}

\noindent Simplifying gives

\begin{align*}
(-1)^{[a]+[b]} q^{(\varepsilon_a,\varepsilon_b)} X^a_b &{A^{(l)}}^b_a
  - q^{(\varepsilon_b, \varepsilon_b)} {A^{(l)}}^b_a X^a_b \notag \\ =
  &\bigl( (-1)^{[a]} - \delta^a_{\overline{b}} q^{(\rho, \varepsilon_a
  - \varepsilon_b)} \bigr) q^{h_{\varepsilon_a}} ({A^{(l)}}^a_a -
  {A^{(l)}}^b_b) \notag \\ & \qquad +(q-q^{-1})\sum_{\varepsilon_a<
  \varepsilon_c < \varepsilon_b} \bigl( (-1)^{[c]} -
  \delta^b_{\overline{c}} q^{(\rho, \varepsilon_c -\varepsilon_b)}
  \bigr) {A^{(l)}}^c_a X^a_c \notag \\ & \qquad +(q-q^{-1})
  \delta^a_{\overline{b}} \sum_{\varepsilon_a< \varepsilon_c <
  \varepsilon_b} (-1)^{[b][c]} \xi_b \xi_c q^{(\rho, \varepsilon_c -
  \varepsilon_b)} X^a_c {A^{(l)}}^{\overline{a}}_{\overline{c}}.
\end{align*}

\noindent Remembering that $\varepsilon_a < \varepsilon_b$, we apply
this to the highest weight state $|\Lambda \rangle$ to obtain

\begin{multline} \label{AX}
-q^{(\varepsilon_b,\varepsilon_b)} {A^{(l)}}^b_a X^a_b |\Lambda
\rangle = q^{(\Lambda, \varepsilon_a)} \bigl( (-1)^{[a]} -
\delta^a_{\overline{b}} q^{2(\rho, \varepsilon_a)} \bigr) (t_a^{(l)} -
t_b^{(l)}) |\Lambda \rangle\\ + (q-q^{-1})\sum_{\varepsilon_a <
  \varepsilon_c < \varepsilon_b} \bigl( (-1)^{[c]} -
\delta^b_{\overline{c}} q^{2(\rho, \varepsilon_c)}\bigr) {A^{(l)}}^c_a
X^a_c |\Lambda \rangle.
\end{multline}


The next step is to calculate ${A^{(l)}}^b_a X^a_b |\Lambda
\rangle$ for $\varepsilon_a < \varepsilon_b$.  It is first convenient
to order the indices according to $b>c \Leftrightarrow \varepsilon_b <
\varepsilon_c$.  With this ordering we say an element $a >0$ if
$\varepsilon_a <0$, $a=0$ if $\varepsilon_a =0$, and $a<0$ if
$\varepsilon_a >0$.  Using this convention, it is apparent the
solution to \eqref{AX} will be of the form

\begin{equation} \label{soln}
{A^{(l)}}^b_a X^a_b |\Lambda \rangle = q^{(\Lambda, \varepsilon_a)}
   (-1)^{[a]} \sum_{a>c \geq b} \alpha^a_{bc} (t_a^{(l)} - t_c^{(l)})
   |\Lambda \rangle,
\end{equation}

\noindent where $\alpha^a_{bc}$ is a function of $a,b$ and $c$. Now
from equation \eqref{AX} we have

\begin{align*}
(q-q^{-1}) \sum_{a>c>b}&(-1)^{[c]}{A^{(l)}}^c_a X^a_c |\Lambda \rangle
\notag\\ &= -q^{(\varepsilon_b,\varepsilon_b)} {A^{(l)}}^b_a X^a_b
|\Lambda \rangle + (q-q^{-1}) \sum_{a>c>b} \delta^b_{\overline{c}}
q^{-2(\rho, \varepsilon_b)} {A^{(l)}}^c_a X^a_c |\Lambda \rangle
\notag \\ & \qquad - (-1)^{[a]} q^{(\Lambda, \varepsilon_a)} \bigl( 1
- \delta^a_ {\overline{b}} (-1)^{[a]} q^{2(\rho, \varepsilon_a)}
\bigr) (t_a^{(l)} - t_b^{(l)}) |\Lambda \rangle \notag \\ &=
-q^{(\varepsilon_{b+1},\varepsilon_{b+1})} {A^{(l)}}^{b+1}_a X^a_{b+1}
|\Lambda \rangle \notag \\ & \qquad + (q-q^{-1}) \sum_{a>c>b+1}
\delta^{b+1}_{\overline{c}} q^{-2(\rho, \varepsilon_{b+1})}
{A^{(l)}}^c_a X^a_c |\Lambda \rangle \notag\\ & \qquad - (-1)^{[a]}
q^{(\Lambda, \varepsilon_a)} \bigl( 1 - \delta^{b+1}_{ \overline{a}}
(-1)^{[a]} q^{2(\rho, \varepsilon_a)} \bigr) (t_a^{(l)} -
t_{b+1}^{(l)}) |\Lambda \rangle \notag \\ & \qquad + (q-q^{-1})
(-1)^{[b+1]} {A^{(l)}}^{b+1}_a X^a_{b+1}|\Lambda \rangle.
\end{align*} 

\noindent Substituting in the form of the solution given in equation
\eqref{soln} produces

\begin{align}
q^{(\varepsilon_b,\varepsilon_b)} &\sum_{a>d\geq b} \alpha^a_{bd}
  (t_a^{(l)} - t_d^{(l)}) |\Lambda \rangle \notag \\ = \bigl(
  &q^{(\varepsilon_{b+1},\varepsilon_{b+1})} - (q-q^{-1}) (-1)^{[b+1]}
  \bigr) \sum_{a>d \geq b+1} \alpha^a_{(b+1)d} (t_a^{(l)} -
  t_d^{(l)}) |\Lambda \rangle \notag \\ 
  -&\bigl(1 -\delta^a_{\overline{b}} (-1)^{[a]} q^{2(\rho,\varepsilon_a)} 
  \bigr) (t_a^{(l)} - t_b^{(l)}) |\Lambda \rangle + \bigl(1 - \delta^a_
  {\overline{b+1}} (-1)^{[a]} q^{2(\rho,\varepsilon_a)} \bigr)
  (t_a^{(l)} - t_{b+1}^{(l)}) |\Lambda \rangle \notag \\ +& (q-q^{-1})
  \sum_{a>c>b} \delta^b_{\overline{c}} q^{-2(\rho, \varepsilon_b)}
  \sum_{a>d \geq c} \alpha^a_{\overline{b}\,d} (t_a^{(l)} - t_d^{(l)})
  |\Lambda \rangle \notag \\ -& (q-q^{-1}) \sum_{a>c>b+1}
  \delta^{b+1}_{\overline{c}} q^{-2(\rho, \varepsilon_{b+1})}
  \sum_{a>d \geq c} \alpha^a_{(\overline{b+1})d} (t_a^{(l)} -
  t_d^{(l)}) |\Lambda \rangle . \label{qaa}
\end{align}

\

\noindent Set

\begin{equation*}
\alpha^a_{bd} = \bar{\alpha}_{bd} (1-\delta^a_{\overline{d}}
  (-1)^{[a]} q^{2(\rho, \varepsilon_a)}).
\end{equation*}

\noindent Then from equation \eqref{qaa} we obtain

\begin{equation*}
\bar{\alpha}_{bb} = - q^{-(\varepsilon_b,\varepsilon_b)}
\end{equation*}

\noindent and

\begin{align*}
\bar{\alpha}_{b(b+1)} &= q^{-(\varepsilon_b,\varepsilon_b)} \Bigl[
  \bigl( q^{(\varepsilon_{b+1},\varepsilon_{b+1})} - (q-q^{-1})
  (-1)^{[b+1]} \bigr) \bar{\alpha}_{(b+1)(b+1)} +1 \\ &\hspace{7cm} +
  (q-q^{-1}) \delta^b_{\overline{b+1}} q^{-2(\rho, \varepsilon_b)}
  \bar{\alpha}_{\overline{b}(b+1)} \Bigr] \\ &=
  q^{-(\varepsilon_b,\varepsilon_b)-(\varepsilon_{b+1},\varepsilon_{b+1})}
  (q-q^{-1}) \bigl( (-1)^{[b+1]} - \delta^b_{\overline{b+1}}
  q^{-2(\rho, \varepsilon_b)} \bigr).
\end{align*}

\

\noindent To simplify this expression note that $q^{2(\rho,
\varepsilon_{b+1} - \varepsilon_b)} =
q^{-(\varepsilon_b,\varepsilon_b) - (\varepsilon_{b+1},
\varepsilon_{b+1})}$ in all cases except for $[b]=0, \,b=l,\, m=2l$,
in which case $q^{2(\rho, \varepsilon_{b+1} - \varepsilon_b)} = q^2
q^{-(\varepsilon_b,\varepsilon_b) - (\varepsilon_{b+1},
\varepsilon_{b+1})}$.  However $[b]=0,\,b=l,\, m=2l$ if and only if
$\delta^b_{\overline{b+1}}=1$, and in that case we find
$\bar{\alpha}_{b(b+1)}=0$.  Hence for all values of $b$ we can write

\begin{equation*}
\bar{\alpha}_{b (b+1)}= (q-q^{-1}) q^{-2(\rho,\varepsilon_b)} \bigl(
  (-1)^{[b+1]}q^{2(\rho,\varepsilon_{b+1})} -
  \delta^b_{\overline{b+1}} \bigr).
\end{equation*}

Now that we have found $\bar{\alpha}_{bb}$ and
$\bar{\alpha}_{b(b+1)}$, they can be used to calculate the remaining
$\bar{\alpha}_{bd}$.  From equation \eqref{qaa} we observe that if
$d>b+1$ then

\begin{align} 
\bar{\alpha}_{bd} &= q^{-(\varepsilon_b,\varepsilon_b)} \bigl
  (q^{(\varepsilon_{b+1}, \varepsilon_{b+1})} - (q-q^{-1})
  (-1)^{[b+1]} \bigr) \bar{\alpha}_{(b+1)d} \nonumber\\ 
& \qquad + (q-q^{-1}) q^{-(\varepsilon_b, \varepsilon_b)} \sum_{d \geq c >b}
  \delta^b_{\overline{c}} q^{-2(\rho, \varepsilon_b)}
  \bar{\alpha}_{\overline {b} d} \nonumber\\ 
& \qquad - (q-q^{-1}) q^{-(\varepsilon_b, \varepsilon_{b})} 
  \sum_{d \geq c >b+1}
  \delta^{b+1}_{\overline{c}} q^{-2(\rho, \varepsilon_{b+1})}
  \bar{\alpha}_{(\overline{b+1}) d} \label{bd1}.
\end{align}

\

\noindent Now define $\theta_{xy}$ by

\begin{equation*}
\theta_{xy} =
\begin{cases} 1 \qquad &x<y,\\
0 & x \geq y.
\end{cases}
\end{equation*}

\noindent Then equation \eqref{bd1} can be rewritten as

\begin{multline} \label{bad}
\bar{\alpha}_{bd} = q^{-(\varepsilon_b,\varepsilon_b)} \bigl(
  q^{(\varepsilon _{b+1}, \varepsilon_{b+1})} - (q-q^{-1})
  (-1)^{[b+1]} \bigr) \bar{\alpha}_ {(b+1)d} \\ + (q-q^{-1})
  q^{-(\varepsilon_b, \varepsilon_b)} q^{2(\rho, \varepsilon_c)}
  \bigl( \theta_{bc} \theta_{c(d+1)} \delta^b_ {\overline{c}} -
  \theta_ {(b+1)c} \theta_{c(d+1)} \delta^{b+1}_{\overline{c}} \bigr)
  \bar{\alpha}_ {cd}, \quad d > a+1.
\end{multline}

\noindent Consider $\bar{\alpha}_{bd}$ for any $b >l$.  Both
$\theta_{b\overline{b}}$ and $\theta_{(b+1)(\overline{b+1})}$ will equal
$0$, so

\begin{align*}
\bar{\alpha}_{bd} &= q^{-(\varepsilon_b,\varepsilon_b)}
  \bigl(q^{(\varepsilon _{b+1}, \varepsilon_{b+1})} - (q-q^{-1})
  (-1)^{[b+1]} \bigr) \bar{\alpha}^a_{(b+1)d} \\ &=
  q^{-(\varepsilon_b, \varepsilon_b)} q^{-(\varepsilon_{b+1},
  \varepsilon_ {b+1})} \bar{\alpha}_{(b+1)d} \\ &= q^{2(\rho,
  \varepsilon_{b+1} - \varepsilon_b)} \bar{\alpha}_{(b+1)d}.
\end{align*}

\noindent Since

\begin{equation*}
\bar{\alpha}_{(d-1)d} =(-1)^{[d]} (q-q^{-1}) q^{2(\rho, \varepsilon_d-
  \varepsilon_{d-1})},
\end{equation*}

\noindent we obtain

\vspace{-1cm}

\begin{equation*} 
\bar{\alpha}_{bd} =(-1)^{[d]} (q-q^{-1}) q^{2(\rho, \varepsilon_{d}-
  \varepsilon_b)}, \quad d>b>l.
\end{equation*}

\noindent Substituting this together with our expression for
$\bar{\alpha}_{bb}$ into equation \eqref{bad}, we find

\begin{align}
\bar{\alpha}_{bd} = &q^{-(\varepsilon_b,\varepsilon_b)} \bigl(
  q^{-(\varepsilon_{b+1}, \varepsilon_{b+1})} -
  \delta^{b+1}_{\overline{b+1}} (q-q^{-1}) \bigr) \bar{\alpha}_{(b+1)d}\notag\\
+& (q-q^{-1})^2 q^{-(\varepsilon_b, \varepsilon_b)}
  (-1)^{[d]} q^{2(\rho, \varepsilon_d)} \bigl( \theta_{b\overline{b}}
  \theta_{\overline {b}d} - \theta_{(b+1)(\overline{b+1})}
  \theta_{(\overline{b+1})d} \bigr)\notag\\ -& (q-q^{-1})
  q^{-(\varepsilon_b, \varepsilon_b)} q^{-(\varepsilon_d,
  \varepsilon_d)} q^{2(\rho, \varepsilon_d)} ( \delta^{\overline{b}}_d
  - \delta^{\overline{b+1}}_d), \; d > b+1. \label{um}
\end{align}

\noindent But for $d>b+1$

\begin{align*}
\theta_{b\overline{b}} \theta_{\overline{b}d} -
  \theta_{(b+1)(\overline{b+1})} \theta_{(\overline{b+1})d} &=
  \delta^b_l \theta_{\overline{l}d} - \delta^{\overline{b}}_d
  \theta_{bl} \notag \\ &= \delta^b_l (1-\delta^d_{\overline{l}}) -
  \delta^{\overline{b}}_d (1- \delta^b_l) \notag \\ &= \delta^b_l
  -\delta^{\overline{b}}_d.
\end{align*}

\noindent Also, $-[(-1)^{[d]} (q-q^{-1}) +
q^{-(\varepsilon_d,\varepsilon_d)}] \delta^{\overline{b}}_d = -
q^{(\varepsilon_d,\varepsilon_d)} \delta^{\overline{b}}_d$, so
equation \eqref{um} reduces to

\begin{align*}
\bar{\alpha}_{bd} &= \bigl( q^{2(\rho, \varepsilon_{b+1} -
  \varepsilon_b)} q^{-2\delta^b_{\overline{b+1}}}
  -\delta^{b+1}_{\overline{b+1}} q^{-1} (q-q^{-1}) \bigr)
  \bar{\alpha}_{(b+1)d} + \delta^b_l q^{-1}(q-q^{-1})^2 (-1)^{[d]}
  q^{2(\rho, \varepsilon_d)} \\ & \quad - \delta^{\overline{b}}_d
  (q-q^{-1}) q^{2(\rho, \varepsilon_d)} + \delta^{\overline{b+1}}_d
  (q-q^{-1}) q^{2(\rho, \varepsilon_{b+1}- \varepsilon_{b})}
  q^{-2\delta^b_{\overline{b+1}}} q^{2(\rho, \varepsilon_d)} \\ 
&=  \bigl( q^{2(\rho, \varepsilon_{b+1} - \varepsilon_b)}
  q^{-2\delta^b_{\overline{b+1}}} -\delta^{b+1}_{\overline{b+1}}
  q^{-1} (q-q^{-1}) \bigr) \bar{\alpha}_{(b+1)d} + \delta^b_l
  q^{-1}(q-q^{-1})^2 (-1)^{[d]} q^{2(\rho, \varepsilon_d)} \\ & \quad
  + (q-q^{-1}) q^{-2(\rho, \varepsilon_b)} (\delta^{\overline{b+1}}_d
  - \delta^{\overline{b}}_d) ,\qquad d>b+1 .
\end{align*}

\noindent Recall that for $b > l$ we have

\begin{equation*}
\bar{\alpha}_{bd} = (-1)^{[d]} (q-q^{-1}) q^{2(\rho,\varepsilon_d
-\varepsilon_b)}, \quad d>b.
\end{equation*}

\noindent Then when $b=l$ we find

\begin{align*}
\bar{\alpha}_{bd} &= \bigl(q^{2(\rho, \varepsilon_{b+1} -
  \varepsilon_b)} q^{-2\delta^b_{\overline{b+1}}} -
  \delta^{b+1}_{\overline{b+1}} q^{-1} (q-q^{-1}) \bigr) (-1)^{[d]}
  (q-q^{-1}) q^{2(\rho,\varepsilon_{d}- \varepsilon_{b+1})} \\ & \quad
  + q^{-1}(q-q^{-1})^2 (-1)^{[d]} q^{2(\rho, \varepsilon_d)} -
  (q-q^{-1}) q^{-2(\rho, \varepsilon_b)} \delta^{\overline{l}}_d
  \notag \\ &= (-1)^{[d]} (q-q^{-1}) q^{2(\rho,\varepsilon_{d}-
  \varepsilon_{b})}\notag\\ &\hspace{2cm} \Bigl[
  \delta^{b+1}_{\overline{b+1}} \bigl( 1 - (q-q^{-1}) + (q-q^{-1})
  \bigr) + \delta^b_{\overline{b+1}} \bigl( q^{-2} + q^{-1} (q-q^{-1})
  \Bigr] \notag \\ & \quad - (q-q^{-1}) q^{-2(\rho, \varepsilon_b)}
  \delta_d^{\overline{l}} \\ &= (q-q^{-1}) q^{-2(\rho,
  \varepsilon_{b})} \bigl( (-1)^{[d]} q^{2(\rho,\varepsilon_{d})} -
  \delta^b_{\overline{d}} \bigr)
\end{align*}

\noindent for all $d > b +1$.  Comparing this with our earlier results
for $d = b+1$ and $b>l$, we have

\begin{equation*}
\bar{\alpha}_{bd} = (q-q^{-1}) q^{-2(\rho, \varepsilon_{b})} \bigl(
  (-1)^{[d]} q^{2(\rho,\varepsilon_{d})} - \delta^b_{\overline{d}}
  \bigr), \quad \forall b \geq l,\; d > b.
\end{equation*}

\noindent But for $b < l$ we know

\begin{align*}
\bar{\alpha}_{bd} &= q^{2(\rho, \varepsilon_{b+1}-\varepsilon_b)}
  \bar{\alpha}_{(b+1)d} + (q-q^{-1}) q^{-2(\rho,\varepsilon_b)}
  (\delta^{\overline{b+1}}_d - \delta^{\overline{b}}_d), \qquad d>b+1
  .
\end{align*}

\noindent Hence for all $b$ we obtain

\begin{align*}
\bar{\alpha}_{bd} &= (q-q^{-1}) q^{-2(\rho, \varepsilon_{b})} \bigl(
  (-1)^{[d]} q^{2(\rho,\varepsilon_{d})} - \sum_{c=b}^{d-1}
  \delta^{\overline {c}}_d + \sum_{c=b}^{d-2}
  \delta^{\overline{c+1}}_d \bigr) \\ &= (q-q^{-1}) q^{-2(\rho,
  \varepsilon_{b})} \bigl( (-1)^{[d]} q^{2(\rho,\varepsilon_{d})} -
  \delta^{\overline{b}}_d \bigr), \qquad d>b.
\end{align*}

\noindent Thus for all $a>b$

\begin{equation}
{A^{(l)}}^b_a X^a_b |\Lambda \rangle = q^{(\Lambda, \varepsilon_a)}
   (-1)^{[a]} \sum_{a>c \geq b} \alpha^a_{bc} (t_a^{(l)} - t_c^{(l)})
   |\Lambda \rangle, \label{al}
\end{equation}

\noindent where $\alpha^a_{bc}$ is given by

\begin{equation*}
\alpha^a_{bc} =
  \begin{cases}
  - q^{-(\varepsilon_b,\varepsilon_b)} (1 -\delta^a_{\overline{b}}
    (-1)^{[a]} q^{2(\rho, \varepsilon_a)}), & c=b, \\ (q-q^{-1})
    q^{-2(\rho, \varepsilon_{b})} \bigl( (-1)^{[c]} q^{2(\rho,
    \varepsilon_{c})} - \delta^{\overline{b}}_c
    \bigr)(1-\delta^a_{\overline{c}} (-1)^{[a]} q^{2(\rho,
    \varepsilon_a)}), \quad & c>b.
  \end{cases}
\end{equation*}

\subsection{Constructing the Perelomov-Popov matrix equation}

\noindent The expression (\ref{al}) can now be substituted into the equation

\begin{equation*}
t^{(l)}_a |\Lambda \rangle = t^{(l-1)}_a t^{(1)}_a |\Lambda \rangle +
  \sum_{\varepsilon_a<\varepsilon_b} (-1)^{[a]+[b]}
  q^{(\Lambda,\varepsilon_a)} {A^{(l-1)}}^b_a X^a_b |\Lambda \rangle
\end{equation*}

\noindent to find a matrix equation for the various $t_a^{(l)}$.  The
matrix factor is an analogue of the Perelomov-Popov matrix introduced
in \cite{Perelomov} and \cite{Popov}, which was used to calculate
the eigenvalues of the Casimir invariants of various classical Lie
algebras.

First recall that

\begin{equation*}
A^b_a = (1+\delta^a_b) X^b_a + (q-q^{-1}) \sum_{c \geq a,b}
(-1)^{([a]+[c])([b]+[c])} X^c_a X^b_c,
\end{equation*}

\noindent where

\begin{equation*}
X^b_a =
\begin{cases}
\frac{q^{h_{\varepsilon_a}} -I}{q-q^{-1}}, &a=b, \\ (-1)^{[b]}
q^{h_{\varepsilon_a}} \hat{\sigma}_{ba}, \quad &\varepsilon_a <
\varepsilon_b, \\ (-1)^{[b]} \hat{\sigma}_{ba} q^{h_{\varepsilon_b}},
&\varepsilon_a > \varepsilon_b.
\end{cases}
\end{equation*}

\noindent Then

\begin{alignat*}{2}
&&\;A^a_a |\Lambda \rangle&= 2 X^a_a |\Lambda \rangle + (q-q^{-1})
X^a_a X^a_a |\Lambda \rangle \notag \\ &&&= (q-q^{-1})^{-1}
(2(q^{h_{\varepsilon_a}}-1) + (q^{h_{\varepsilon_a}}-1)^2)
|\Lambda\rangle. \\ &\therefore & t^{(1)}_a & =\frac{q^{2(\Lambda,
    \varepsilon_a)}-1}{q-q^{-1}}.
\end{alignat*}

\noindent Hence we obtain

\begin{align*}
t^{(l)}_a &= \frac{(q^{2(\Lambda, \varepsilon_a)} -1)}{(q-q^{-1})}
  t^{(l-1)}_a \notag \\ &\qquad+ \sum_{b<a} (-1)^{[a]+[b]}
  q^{(\Lambda,\varepsilon_a)} \bigl(q^ {(\Lambda,\varepsilon_a)}
  (-1)^{[a]} \sum_{b \leq c < a} \alpha^a_{bc} (t^{(l-1)}_a -
  t^{(l-1)}_c) \bigr) \notag \\ &= \frac{(q^{2(\Lambda,
  \varepsilon_a)} -1)}{(q-q^{-1})} t^{(l-1)}_a \notag\\ &\qquad
  -q^{2(\Lambda,\varepsilon_a)} \sum_{b<a} (-1)^{[b]}
  q^{-(\varepsilon_b, \varepsilon_b)}
  (1-\delta^a_{\overline{b}}(-1)^{[a]} q^{2(\rho, \varepsilon_a)})
  (t^{(l-1)}_a - t^{(l-1)}_b) \notag \\ & \qquad + (q-q^{-1})
  q^{2(\Lambda,\varepsilon_a)} \sum_{c < b<a} (-1)^{[c]}
  q^{-2(\rho,\varepsilon_{c})} (1-\delta^a_{\overline{b}} (-1)^{[a]}
  q^{2(\rho, \varepsilon_a)}) \notag \\ & \hspace{6cm} ((-1)^{[b]}
  q^{2(\rho,\varepsilon_{b})}-\delta^b_{\overline{c}}) (t^{(l-1)}_a -
  t^{(l-1)}_b) .
\end{align*}

Now consider the function $\gamma_b$ defined by:

\begin{equation*}
\gamma_b = (-1)^{[b]} q^{-(\varepsilon_b,\varepsilon_b)} - (q-q^{-1})
  \sum_{c < b} (-1)^{[c]} q^{-2(\rho, \varepsilon_{c})} \bigl(
  (-1)^{[b]} q^{2(\rho,\varepsilon_{b})} - \delta^b_{\overline{c}}
  )\bigr).
\end{equation*}

\noindent We evaluate this for all $b$, remembering that $C(\Lambda_0) =
(\delta_1, \delta_1+2\rho) = m-n-1$ and 

\begin{equation*}
\rho = \frac{1}{2} \sum_{i=1}^l (m-2i) \varepsilon_i + \frac{1}{2}
\sum_{\mu=1}^k (n-m+2-2\mu) \delta_\mu.
\end{equation*}

\noindent We find

\begin{equation*}
\gamma_b = (-1)^{[b]} q^{2(\rho, \varepsilon_b)} q^{-C(\Lambda_0)}
\end{equation*}

\noindent for all values of $b$.  We also consider the function

\begin{equation*}
\beta_a = 1 - (q-q^{-1}) \sum_{b<a} \gamma_b \bigl( 1 -
  \delta^a_{\overline{b}} (-1)^{[a]} q^{2(\rho, \varepsilon_a)}\bigr),
\end{equation*}

\noindent so that

\begin{equation} \label{tal}
t_a^{(l)} = \frac{(q^{2(\Lambda, \varepsilon_a)} \beta_a
  -1)}{(q-q^{-1})} t_a^{(l-1)} + q^{2(\Lambda, \varepsilon_a)}
  \sum_{b<a} \gamma_b \bigl( 1- \delta^a_{\overline{b}} (-1)^{[a]}
  q^{2(\rho, \varepsilon_a)}\bigr) t_b^{(l-1)}.
\end{equation}

\noindent Again, by considering the various cases individually we find

\begin{equation*}
\beta_a= q^{(\varepsilon_a, 2\rho +\varepsilon_a) - C(\Lambda_0)}
\end{equation*}

\noindent for any $a$, regardless of whether $m$ is even or odd.  Substituting
this result together with that for $\gamma_b$ into equation
(\ref{tal}) gives

\begin{multline*}
t_a^{(l)} = \frac{(q^{(\varepsilon_a,2\Lambda+ 2\rho +\varepsilon_a) -
  C(\Lambda_0)}-1)}{(q-q^{-1})} t_a^{(l-1)} \\ + q^{(2\Lambda,
  \varepsilon_a)-C(\Lambda_0)} \sum_{b<a} (-1)^{[b]} q^{(2\rho,
  \varepsilon_b)} \bigl( 1-\delta^a_{\overline{b}} (-1)^{[a]}
  q^{(2\rho, \varepsilon_a)}\bigr) t_b^{(l-1)}.
\end{multline*}

\noindent This can be written in the matrix form

\begin{equation*}
\underline{t}^{(l)} = M \underline{t}^{(l-1)},
\end{equation*}

\noindent where $M$ is a lower triangular matrix with entries

\begin{equation*}
M_{ab} =
\begin{cases}
0, &a<b \\ (q-q^{-1})^{-1}(q^{(\varepsilon_a,2\Lambda+ 2\rho
+\varepsilon_a) - C(\Lambda_0)}-1), &a=b, \\ q^{(2\Lambda,
\varepsilon_a)-C(\Lambda_0)} \bigl( (-1)^{[b]} q^{(2\rho,
\varepsilon_b)} -\delta^a_{\overline{b}} \bigr), &a>b.
\end{cases}
\end{equation*}

\noindent Then we have

\begin{equation*}
\underline{t}^{(l)} = M^{l} \underline{t}^{(0)}, \qquad \text{with}
  \quad t^{(0)}_a =1 \quad \forall a,
\end{equation*}

\noindent where $M$ is an analogue of the Perelomov-Popov matrix.

\subsection{Solving the matrix equation}

\noindent This matrix equation for $t_a^{(l)}$ can now be used to
calculate the eigenvalues of $C_l$.  Loosely speaking, the problem
reduces to diagonalising the matrix $M$.
Recall

\begin{equation*}
C_l = \sum_a (-1)^{[a]} q^{(2\rho, \varepsilon_a)} {A^{(l)}}^a_a.
\end{equation*}

\noindent Denote the eigenvalue of $C_l$ on $V(\Lambda)$ as
$\chi_{\Lambda}(C_l)$. Then we have

\begin{equation*}
\chi_{\Lambda}(C_l) = \sum_a (-1)^{[a]} q^{(2\rho, \varepsilon_a)}
  t_a^{(l)} = \sum_{a,b} (-1)^{[a]} q^{(2\rho, \varepsilon_a)}
  (M^l)_{ab}.
\end{equation*}

\noindent To calculate this we wish to diagonalise $M$.  We assume the
eigenvalues of $M$,

\begin{equation*}
\alpha_a^{\Lambda} = \frac{(q^{(\varepsilon_a,2\Lambda+ 2\rho
  +\varepsilon_a) - C(\Lambda_0)}-1)}{(q-q^{-1})},
\end{equation*}

\noindent are distinct.  Then we need a matrix $N$ satisfying

\begin{equation*}
(N^{-1}MN)_{ab} = \delta^a_b \alpha_a^{\Lambda},
\end{equation*}

\noindent which implies

\begin{equation} \label{chi}
\chi_{\Lambda}(C_l) = \sum_{a,b,c} (-1)^{[a]} q^{(2\rho,
  \varepsilon_a)} ( \alpha_b^{\Lambda} )^l N_{ab} (N^{-1})_{bc}.
\end{equation}

Now

\begin{equation*}
(MN)_{ab} = \alpha_b^{\Lambda} N_{ab}.
\end{equation*}

\noindent Substituting in the values for $M_{ab}$ gives

\begin{equation} \label{alphaN}
\alpha_a^{\Lambda} N_{ab} + q^{(2\Lambda, \varepsilon_a)-C(\Lambda_0)}
  \sum_{c<a} \bigl( (-1)^{[c]} q^{(2\rho, \varepsilon_c)} -
  \delta^a_{\overline{c}} \bigr) N_{cb} = \alpha_b^{\Lambda} N_{ab}.
\end{equation}

\noindent Since the eigenvalues $\alpha_a^{\Lambda}$ are distinct,
this implies

\begin{equation*}
N_{ab} = 0, \qquad \forall a < b.
\end{equation*}

\noindent Set

\begin{equation} \label{P}
P_{ab} = \sum_{c\leq a} (-1)^{[c]} q^{(2\rho, \varepsilon_c)} N_{cb}.
\end{equation}

\noindent Then equation (\ref{alphaN}) becomes

\begin{alignat*}{2}
&&&(\alpha_b^{\Lambda}- \alpha_a^{\Lambda}) N_{ab} = q^{(2\Lambda,
  \varepsilon_a) -C(\Lambda_0)} P_{(a-1)\,b} - \theta_{0a} q^{(2\Lambda,
  \varepsilon_a)-C(\Lambda_0)} N_{\overline{a}b} \notag \\
  &\Rightarrow \qquad && (\alpha_b^{\Lambda} -
  \alpha_a^{\Lambda})(-1)^{[a]} q^{(-2\rho,\varepsilon_a)}
  (P_{ab}-P_{(a-1)\,b}) \notag \\ 
&&&\hspace{35mm}= q^{(2\Lambda,
  \varepsilon_a)-C(\Lambda_0)} P_{a-1\,b} - \theta_{0a} q^{(2\Lambda,
  \varepsilon_a)-C(\Lambda_0)} N_{\overline{a}b},
\end{alignat*}

\noindent which simplifies to

\begin{equation*} 
P_{ab} = \frac{(\alpha_b^{\Lambda} - \alpha_a^{\Lambda} + (-1)^{[a]}
  q^{2(\Lambda+\rho, \varepsilon_a) -
  C(\Lambda_0)})}{(\alpha_b^{\Lambda}- \alpha_a^{\Lambda})} P_{(a-1)\,b}
  - \frac{\theta_{0a}(-1)^{[a]} q^{2(\Lambda + \rho, \varepsilon_a)
  -C(\Lambda_0)}}{(\alpha_b^{\Lambda} - \alpha_a^{\Lambda})}
  N_{\overline{a}b}.
\end{equation*}
\

\noindent Set

\begin{equation*}
\psi^b_{a} = \alpha_b^{\Lambda} - \alpha_a^{\Lambda} + (-1)^{[a]}
  q^{2(\Lambda+\rho, \varepsilon_a) - C(\Lambda_0)},
\end{equation*}

\noindent so this becomes

\begin{equation} \label{Pab}
P_{ab} = \frac{\psi^b_a}{(\alpha_b^{\Lambda}- \alpha_a^{\Lambda})}
P_{(a-1)\,b} - \frac{\theta_{0a} (-1)^{[a]} q^{2(\Lambda+\rho,
    \varepsilon_a) - C(\Lambda_0)}} {(\alpha_b^{\Lambda} -
  \alpha_a^{\Lambda})} N_{\overline{a}b}.
\end{equation}

\noindent Without loss of generality we can choose $N_{aa} = 1\;
\forall a$, so $P_{bb} = (-1)^{[b]} q^{2(\rho, \varepsilon_b)}$.  Then
in the cases $0 \geq a > b$ and $a > b \geq 0$ the last term in
equation (\ref{Pab}) vanishes, giving

\begin{equation*}
P_{ab} = (-1)^{[b]} q^{2(\rho, \varepsilon_b)} \prod_{c=b+1}^a
  \frac{\psi^b_c} {(\alpha_b^{\Lambda}- \alpha_c^{\Lambda})}.
\end{equation*}

\noindent Similarly, for $a > \overline{b} > 0$ we obtain

\begin{equation} \label{abbar0}
P_{ab} = P_{\overline{b}b} \prod_{c=\overline{b}+1}^a \frac{\psi^b_c}
  {(\alpha_b^{\Lambda}- \alpha_c^{\Lambda})} .
\end{equation}

\noindent It remains to find $P_{ab}$ for $\overline{b} \geq a >0$.
In this case, the last term in equation (\ref{Pab}) contributes,
giving

\begin{multline} \label{P_ab}
P_{ab} = (-1)^{[b]} q^{2(\rho, \varepsilon_b)} \prod_{c=b+1}^a
  \frac{\psi^b_c} {(\alpha_b^{\Lambda}- \alpha_c^{\Lambda})} -
  \frac{(-1)^{[a]} q^{2(\Lambda+ \rho, \varepsilon_a)
  -C(\Lambda_0)}}{(\alpha_b^{\Lambda}-\alpha_a^{\Lambda})}
  N_{\overline{a}b} \\ - \sum_{d=\overline{l}}^{a-1}\frac{(-1)^{[d]}
  q^{2(\Lambda +\rho,\varepsilon_d)
  -C(\Lambda_0)}}{(\alpha_b^{\Lambda} - \alpha_d^{\Lambda})}
  N_{\overline{d}b} \prod_{c=d+1}^a \frac{\psi^b_c}
  {(\alpha_b^{\Lambda}- \alpha_c^{\Lambda})}.
\end{multline}

\noindent Recall that if $b<a<0$, then

\begin{align*}
N_{ab} &=
  \frac{q^{(2\Lambda,\varepsilon_a)-C(\Lambda_0)}}{(\alpha_b^{\Lambda}-
  \alpha_a^{\Lambda})} P_{(a-1)\, b} \notag \\ &= \frac{(-1)^{[b]}
  q^{2(\Lambda,\varepsilon_a)+2(\rho,\varepsilon_b)-C( \Lambda_0)}}
  {(\alpha_b^{\Lambda}- \alpha_a^{\Lambda})} \prod_{c=b+1}^{a-1}
  \frac{\psi^b_c} {(\alpha_b^{\Lambda}- \alpha_c^{\Lambda})}.
\end{align*} 

\noindent Substituting this into equation \eqref{P_ab}, we find

\begin{align*}
P_{\overline{b}b} &= (-1)^{[b]} q^{2(\rho,\varepsilon_b)}
  \prod_{c=b+1}^ {\overline{b}} \frac{\psi^b_c} {(\alpha_b^{\Lambda}-
  \alpha_c^{\Lambda})} - \frac{(-1)^{[b]} q^{-2(\Lambda+\rho,
  \varepsilon_b) -C(\Lambda_0)}} {(\alpha_b^{\Lambda} -
  \alpha_{\overline{b}}^{\Lambda})} \notag \\ &\qquad
  -\sum_{d=\overline{l}}^{\overline{b}-1} \frac{(-1)^{[d]+[b]}
  q^{2(\rho, \varepsilon_d + \varepsilon_b) - 2C(\Lambda_0)}}
  {(\alpha_b^{\Lambda} - \alpha_d^{\Lambda}) (\alpha_b^{\Lambda}
  -\alpha_{\overline{d}}^{\Lambda})} \prod_{c=b+1}^{\overline{d}-1}
  \frac{\psi^b_c} {(\alpha_b^{\Lambda}- \alpha_c^{\Lambda})}
  \prod_{c=d+1}^{\overline{b}}\frac{\psi^b_c} {(\alpha_b^{\Lambda} -
  \alpha_c^{\Lambda})},
\end{align*}

\noindent which can also be simplified to

\begin{equation}
P_{\overline{b}b} \prod_{c=b+1}^{\overline{b}}
  \frac{(\alpha_b^{\Lambda}- \alpha_c^{\Lambda})} {\psi^b_c} =
  (-1)^{[b]} q^{2(\rho,\varepsilon_b)} \Biggl[
  1-\sum_{d=\overline{l}}^{\overline{b}} \frac{(-1)^{[d]}
  q^{2(\rho,\varepsilon_d)-2C(\Lambda_0)}} {\psi^b_d
  \psi^b_{\overline{d}}} \prod_{c=\overline{d}+1}^{d-1}
  \frac{(\alpha_b^{\Lambda}-\alpha_c^{\Lambda})}
  {\psi^b_c}\Biggr]\label{bbbar}.
\end{equation}

From this point we will consider the case $m=2l+1$.  This is
marginally more complicated than the case with even $m$. Define $\Phi^b_d$ 
to be

\begin{align*}
\Phi^b_d &= \prod_{c=\overline{l}}^{d-1} \frac{(\alpha_b
  -\alpha_c)(\alpha_b - \alpha_{\overline{c}})}{\psi^b_c
  \psi^b_{\overline{c}}} \\ &= \frac{(\alpha_b -
  \alpha_{d-1})(\alpha_b - \alpha_{\overline{d-1}})} {\psi^b_{d-1}
  \psi^b_{\overline{d-1}}} \Phi^b_{d-1}, \quad \Phi^b_
  {\overline{l}}=1.
\end{align*}

\noindent Then $P_{\overline{b}b}$ can be written as

\begin{equation*}
P_{\overline{b}b} = (-1)^{[b]} q^{2(\rho, \varepsilon_b)}
  \prod_{^{c=b+1}_{c \neq 0}}^{\overline{b}}
  \frac{\psi^b_c}{(\alpha_b^{\Lambda}-\alpha_c^ {\Lambda})} \Biggl[
  \frac{\psi^b_0}{\alpha_b - \alpha_0} - \sum_{d=\overline
  {l}}^{\overline{b}} \frac{(-1)^{[d]}
  q^{2(\rho,\varepsilon_d)-2C(\Lambda_0)}} {\psi^b_d
  \psi^b_{\overline{d}}} \Phi^b_d \Biggr].
\end{equation*}

\noindent Note that for $c \neq 0$,

\begin{align*}
\psi^b_c &= \frac{q^{-C(\Lambda_0)}}{(q-q^{-1})} \bigl(
  q^{(\varepsilon_b, 2\rho+2\Lambda +\varepsilon_b)} -
  q^{(\varepsilon_c, 2\rho+2\Lambda + \varepsilon_c)}
  +(q-q^{-1})(-1)^{[c]} q^{(\varepsilon_c, 2\rho + 2\Lambda)} \bigr)
  \notag \\ &= \frac{ q^{-C(\Lambda_0)}\,\tilde{\psi}^b_c}{(q-q^{-1})}
\end{align*}

\noindent where

\begin{equation*}
\tilde{\psi}^b_c = q^{(\varepsilon_b, 2\rho +2\Lambda +
\varepsilon_b)} - q^{(\varepsilon_c, 2\rho+2\Lambda-\varepsilon_c)}.
\end{equation*}

\noindent So

\begin{align}
\sum_{d=\overline{l}}^{\overline{b}}
\frac{(-1)^{[d]}q^{2(\rho,\varepsilon_d) -2C(\Lambda_0)}} {\psi^b_d
  \psi^b_{\overline{d}}} \Phi^b_d &=(q-q^{-1})
\sum_{d=\overline{l}}^{\overline{b}} \frac{(-1)^{[d]} (q-q^{-1})
  q^{2(\rho,\varepsilon_d)}} {\tilde{\psi}^b_d
  \tilde{\psi}^b_{\overline{d}}} \Phi^b_d \notag \\
&=(q-q^{-1})\sum_{d=\overline{l}}^{\overline{b}}
\frac{(q^{2(\varepsilon_d, \varepsilon_d)}-1) q^{2(\rho,\varepsilon_d)
    - (\varepsilon_d, \varepsilon_d)}}{\tilde{\psi}^b_d
  \tilde{\psi}^b_{\overline{d}}} \Phi^b_d
\label{sumjlb}
\end{align}

\noindent and

\begin{align*}
\Phi^b_{d+1} &= \frac{(\alpha_b -
  \alpha_d)(\alpha_b-\alpha_{\overline{d}})} {\psi^b_d
  \psi^b_{\overline{d}}} \Phi^b_{d} \\ &= \frac{(q^{(\varepsilon_b,
  \varepsilon_b + 2\rho + 2\Lambda)} - q^{(\varepsilon_d,
  \varepsilon_d + 2\rho + 2\Lambda)})(q^{(\varepsilon_b, \varepsilon_b
  + 2\rho + 2\Lambda)} - q^{(\varepsilon_d, \varepsilon_d - 2\rho -
  2\Lambda)})}{\tilde{\psi}^b_d \tilde{\psi}^b_{\overline{d}}}
  \Phi^b_{d}
\end{align*}

\noindent for $d \geq \overline{l}$. Now

\begin{align*}
&(q^{(\varepsilon_b, \varepsilon_b + 2\rho + 2\Lambda)} -
  q^{(\varepsilon_d, \varepsilon_d + 2\rho +
  2\Lambda)})(q^{(\varepsilon_b, \varepsilon_b + 2\rho + 2\Lambda)} -
  q^{(\varepsilon_d, \varepsilon_d - 2\rho - 2\Lambda)})\notag\\ &
  \hspace{1cm} = q^{2(\varepsilon_d, \varepsilon_d)}
  (q^{(\varepsilon_b, \varepsilon_b + 2\rho + 2\Lambda)} -
  q^{(\varepsilon_d, -\varepsilon_d + 2\rho + 2\Lambda)})
  (q^{(\varepsilon_b, \varepsilon_b + 2\rho + 2\Lambda)} -
  q^{-(\varepsilon_d, \varepsilon_d + 2\rho + 2\Lambda)}) \notag\\ &
  \hspace{1cm} \qquad + q^{2(\varepsilon_b,\varepsilon_b + 2\rho +
  2\Lambda)} (1-q^{2(\varepsilon_d, \varepsilon_d)}) +
  q^{2(\varepsilon_d, \varepsilon_d)}-1 \notag\\
  &\hspace{1cm}=q^{2(\varepsilon_d, \varepsilon_d)} \tilde{\psi}^b_d
  \tilde{\psi} ^b_{\overline{d}} - (q^{2(\varepsilon_b,\varepsilon_b +
  2\rho + 2\Lambda)}-1) (q^{2(\varepsilon_d, \varepsilon_d)}-1).
\end{align*}

\noindent Then, for $d\geq \overline{l}$,

\begin{equation} \label{cancel}
\frac{\Phi^b_{d+1}}{(q^{2(\varepsilon_b,\varepsilon_b + 2\rho
+2\Lambda)}-1)} = \Bigl[ \frac{q^{2(\varepsilon_d,
\varepsilon_d)}}{(q^{2(\varepsilon_b, \varepsilon_b + 2\rho +
2\Lambda)}-1)} - \frac{(q^{2(\varepsilon_d, \varepsilon_d)}-1)}
{\tilde{\psi}^b_d \tilde{\psi}^b_{\overline{d}}} \Bigr] \Phi^b_{d}.
\end{equation}

\noindent Now for $d= \overline{b}$

\begin{align*}
&\frac{(q^{2(\varepsilon_d, \varepsilon_d)}-1)
  q^{2(\rho,\varepsilon_d) -(\varepsilon_d,
    \varepsilon_d)}}{\tilde{\psi}^b_d \tilde{\psi}^b_{\overline {d}}}
\\ & \hspace{25mm} = \frac{(q^{2(\varepsilon_b, \varepsilon_b)}-1)
  q^{2(\rho,\varepsilon_{\overline{b}}) - (\varepsilon_b,
    \varepsilon_b)}} {(q^{(\varepsilon_b, 2\rho+2\Lambda +
    \varepsilon_b)} - q^{-(\varepsilon_b, 2\rho + 2\Lambda +
    \varepsilon_b)}) q^{(\varepsilon_b, 2\rho + 2\Lambda)}
  (q^{(\varepsilon_b, \varepsilon_b)} - q^{-(\varepsilon_b,
    \varepsilon_b)})}\\ & \hspace{25mm} =
\frac{q^{2(\rho,\varepsilon_{\overline{b}}) +(\varepsilon_b,
    \varepsilon_b)}} {(q^{2(\varepsilon_b,\varepsilon_b+2\rho +
    2\Lambda)} - 1)},
\end{align*}

\noindent which can be written as

\begin{equation*}
  \frac{(q^{2(\varepsilon_d, \varepsilon_d)}-1)
  q^{2(\rho,\varepsilon_d) - (\varepsilon_d,
  \varepsilon_d)}}{\tilde{\psi}^b_d \tilde{\psi}^b_{\overline {d}}} =
  \frac{q^{2(\rho,\varepsilon_{\overline{b}-1})-(\varepsilon_
  {\overline{b}-1}, \varepsilon_{\overline{b}-1})}}
  {(q^{2(\varepsilon_b, \varepsilon_b+ 2\rho +2\Lambda) } - 1)}
\end{equation*}

\noindent when $b<l$.  Hence equation \eqref{cancel} can be used to
pairwise cancel the terms in the sum in equation \eqref{sumjlb}.
Adding the first two terms ($d=\overline{b},\, \overline{b}-1$), we
find:

\begin{align*}
q^{2(\rho,\varepsilon_{\overline{b}-1})-(\varepsilon_{\overline{b}-1},
  \varepsilon_{\overline{b}-1})} \Bigl[ \frac{\Phi^b_{\overline{b}}}
  {(q^{2 (\varepsilon_b,\varepsilon_b + 2\rho + 2\Lambda)} - 1)} &+
  \frac{(q^{2 (\varepsilon_{\overline{b}-1},
  \varepsilon_{\overline{b}-1})}-1)} {\tilde {\psi}^b_{\overline{b}-1}
  \tilde{\psi}^b_{b+1}} \Phi^b_{\overline{b}-1} \Bigr] \notag \\ &=
  q^{2(\rho,\varepsilon_{\overline{b}-1})-(\varepsilon_{\overline{b}-1},
  \varepsilon_{\overline{b}-1})}
  \frac{q^{2(\varepsilon_{\overline{b}-1},
  \varepsilon_{\overline{b}-1})}}{(q^{2(\varepsilon_b,\varepsilon_b +
  2\rho + 2\Lambda)}-1)}\Phi^b_{\overline{b}-1} \notag \\ &=
  \frac{q^{2(\rho,\varepsilon_{\overline{b}-2})-(\varepsilon_{\overline{b}-2},
  \varepsilon_{\overline{b}-2})}}{(q^{2(\varepsilon_b,\varepsilon_b +
  2\rho + 2\Lambda)} - 1)}\Phi^b_{\overline{b}-1}.
\end{align*}

\noindent Continuing to apply equation (\ref{cancel}) in this manner
gives

\begin{align}
\sum_{d=\overline{l}}^{\overline{b}} \frac{(q^{2(\varepsilon_d,
  \varepsilon_d)}-1) q^{2(\rho,\varepsilon_d) - (\varepsilon_d,
  \varepsilon_d)}}{\tilde{\psi}^b_d \tilde{\psi}^b_{\overline{d}}}
  \Phi^b_d &=
  \frac{q^{2(\rho,\varepsilon_{\overline{l}})+(\varepsilon_l,
  \varepsilon_l)}}{(q^{2 (\varepsilon_b,\varepsilon_b + \rho +
  \Lambda)} - 1)} \Phi^b_{\overline{l}}\notag \\ &=
  \frac{q^{2l+1-m}}{(q^{2 (\varepsilon_b,\varepsilon_b +2\rho
  +2\Lambda)}-1)}.
  \label{sumprod}
\end{align}

\noindent Hence in the case $m=2l+1$

\begin{equation*}
P_{\overline{b}b} = (-1)^{[b]} q^{2(\rho,\varepsilon_b)}
  \Bigl[\frac{\psi^b_0} {\alpha_b - \alpha_0} - \frac{(q-q^{-1})}
  {(q^{2 (\varepsilon_b,\varepsilon _b + 2\rho + 2\Lambda)}- 1)}
  \Big]\prod_{^{c=b+1}_{c \neq 0}}^{\overline{b}}
  \frac{\psi^b_c}{(\alpha_b^{\Lambda}-\alpha_c^{\Lambda})}.
\end{equation*}

\noindent By substituting in the formulae for $\psi^b_c$ and
$\alpha_b$ and simplifying we obtain

\begin{equation*}
P_{\overline{b}b}= (-1)^{[b]} q^{2(\rho, \varepsilon_b)} \Bigl[1 +
  (q-q^{-1}) \frac{q^{(\varepsilon_b, \varepsilon_b + 2\rho
  +2\Lambda)}}{(q^{2 (\varepsilon_b, \varepsilon_b + 2\rho +
  2\Lambda)}- 1)} \Bigr] \prod_{c=b+1}^{\overline{b}}
  \frac{(q^{(\varepsilon_b, 2\rho +2\Lambda + \varepsilon_b)} -
  q^{(\varepsilon_c, 2\rho+2\Lambda-\varepsilon_c)})}
  {(q^{(\varepsilon_b, 2\rho +2\Lambda
  +\varepsilon_b)}-q^{(\varepsilon_c,
  2\rho+2\Lambda+\varepsilon_c)})},
\end{equation*}

\noindent and thus for $a \geq \overline{b} > 0$

\begin{equation*}
P_{ab} = (-1)^{[b]} q^{2(\rho, \varepsilon_b)} \Bigl[1 + (q-q^{-1})
  \frac {q^{(\varepsilon_b, \varepsilon_b + 2\rho +2\Lambda)}}{(q^{2
  (\varepsilon_b, \varepsilon_b + 2\rho + 2\Lambda)}- 1)} \Bigr]
  \prod_{c=b+1}^{a} \frac{(q^{(\varepsilon_b, 2\rho +2\Lambda +
  \varepsilon_b)} - q^{(\varepsilon_c,
  2\rho+2\Lambda-\varepsilon_c)})} {(q^{(\varepsilon_b, 2\rho +
  2\Lambda+\varepsilon_b)}-q^{(\varepsilon_c,
  2\rho+2\Lambda+\varepsilon_c)})}.
\end{equation*}

\noindent Similarly, we find from equations \eqref{abbar0},
\eqref{bbbar}, \eqref{sumjlb} and \eqref{sumprod} that when $m$ is even
then

\begin{equation*}
P_{ab} = (-1)^{[b]} q^{2(\rho,\varepsilon_b)} \Bigl[1 -
  \frac{q(q-q^{-1})} {(q^{2 (\varepsilon_b,\varepsilon_b + 2\rho +
  2\Lambda)}- 1)} \Big] \prod_ {c=b+1}^{a} \frac{(q^{(\varepsilon_b,
  2\rho +2\Lambda +\varepsilon_b)} - q^{(\varepsilon_c,
  2\rho+2\Lambda-\varepsilon_c)})}{(q^{(\varepsilon_b, 2\rho
  +2\Lambda+\varepsilon_b)}-q^{(\varepsilon_c,2\rho+2\Lambda+\varepsilon_c)})}
\end{equation*}

\noindent for $a \geq \overline{b} > 0$. Hence we have found
expressions for $P_{ab}$ for all $a,b$ satisfying $a \geq \overline{b}
>0$.  At the end of the paper these, together with the earlier results
for $P_{ab}$, will be used to calculate $\chi_\Lambda (C_l)$.


Now we return to the diagonalisation of the matrix $N$.  We know

\begin{equation*}
(N^{-1}M)_{ab} = \alpha_a^\Lambda (N^{-1})_{ab}.
\end{equation*}

\noindent Substituting in the values for $M_{ab}$ gives

\begin{equation} \label{Q} 
\alpha_b^\Lambda (N^{-1})_{ab} + (-1)^{[b]} q^{(2\rho, \varepsilon_b)
  - C(\Lambda_0)} \sum_{c>b} q^{(2\Lambda, \varepsilon_c)} (1-
  \delta^c_{ \overline{b}} (-1)^{[b]} q^{-2(\rho, \varepsilon_b)})
  (N^{-1})_{ac} = \alpha_a^{\Lambda} (N^{-1})_{ab}.
\end{equation}

\noindent Set

\begin{equation*}
\hat{Q}_{ab} = \sum_{c \geq b} q^{2(\Lambda, \varepsilon_c)}
(N^{-1})_{ac}.
\end{equation*}

\noindent We then solve for $\hat{Q}_{ab}$, with the calculations
being very similar to those for $P_{ab}$.  For $0\leq b<a$ and $b < a
\leq 0$ we find

\begin{equation*}
\hat{Q}_{ab} = q^{2(\Lambda, \varepsilon_a)} \prod_{c=b}^{a-1}
  \frac{\psi^a_c}{(\alpha_a^{\Lambda} - \alpha_c^\Lambda)}.
\end{equation*}

\noindent For $m=2l+1$ we obtain

\begin{equation*}
\hat{Q}_{ab} = q^{2(\Lambda,\varepsilon_a)} \Bigl[ 1 + (q-q^{-1})
  \frac {q^{(\varepsilon_a, \varepsilon_a + 2\rho +2\Lambda)}}{(q^{2
  (\varepsilon_a, \varepsilon_a + 2\rho + 2\Lambda)}- 1)} \Bigr]
  \prod_{c=b}^{a-1} \frac{(q^{(\varepsilon_a, 2\rho +2\Lambda +
  \varepsilon_a)} - q^{(\varepsilon_c,
  2\rho+2\Lambda-\varepsilon_c)})} {(q^{(\varepsilon_a, 2\rho
  +2\Lambda+\varepsilon_a)}-q^{(\varepsilon_c, 2\rho
  +2\Lambda+\varepsilon_c)})}.
\end{equation*}

\noindent for $b \leq \overline{a} < 0$.  Similarly, for even $m$ we
find

\begin{equation*}
\hat{Q}_{ab} = q^{2(\Lambda,\varepsilon_a)} \Bigl[ 1 -
  \frac{q(q-q^{-1})} {(q^{2 (\varepsilon_a,\varepsilon_a + 2\rho +
  2\Lambda)}- 1)} \Big] \prod_{c=b}^{a-1} \frac{(q^{(\varepsilon_a,
  2\rho +2\Lambda + \varepsilon_a)} - q^{(\varepsilon_c,
  2\rho+2\Lambda-\varepsilon_c)})} {(q^{(\varepsilon_a, 2\rho
  +2\Lambda+\varepsilon_a)}-q^{(\varepsilon_c, 2\rho
  +2\Lambda+\varepsilon_c)})}.
\end{equation*}

\noindent for $b \leq \overline{a} < 0$.

To use these results to calculate $\chi_\Lambda (C_l)$ we
introduce a new function $Q_{ab}$, defined by:

\begin{equation*}
Q_{ab} = \sum_{c\geq b} (N^{-1})_{ac}.
\end{equation*}

\noindent Then from equations \eqref{chi} and \eqref{P} we deduce

\begin{equation} \label{PQ}
\chi_{\Lambda}(C_l) = \sum_a {(\alpha_a^{\Lambda})}^l
P_{(\mu=\overline{1})\,a} Q_{a\, (\nu = 1)}.
\end{equation}

\noindent However we know

\begin{equation*}
t_a^{(l)} = \frac{q^{2(\Lambda, \varepsilon_a)} - 1}{q-q^{-1}},
\end{equation*}

\noindent and

\begin{align*}
\sum_b (N^{-1})_{ab} t_b^{(1)} &= \sum_{b,c} (N^{-1})_{ab} M_{bc}
  t_c^{(0)} \notag \\ \Rightarrow \quad \sum_b
  (N^{-1})_{ab}\frac{(q^{2(\Lambda,\varepsilon_b)} - 1)} {(q-q^{-1})}
  &= \sum_b (N^{-1}M)_{ab} \notag \\ &= \sum_b \alpha_a^\Lambda
  (N^{-1})_{ab} \notag \\ &= \sum_b (N^{-1})_{ab}
  \frac{(q^{(\varepsilon_a,2\Lambda+2\rho +\varepsilon_a) -
  C(\Lambda_0)}-1)}{(q-q^{-1})}.
\end{align*}

\noindent Thus

\begin{equation*}
Q_{a\, (\nu=1)} = q^{C(\Lambda_0) - (\varepsilon_a, 2\rho + 2\Lambda +
  \varepsilon_a)} \hat{Q}_{a\, (\nu=1)}.
\end{equation*}

\subsection{Explicit formulae for the eigenvalues}

\noindent Substituting our formulae for $ P_{(\mu=\overline{1})\, a}$
and $Q_{a\, (\nu=1)}$ into equation \eqref{PQ}, noting that for $a \neq
0$ exactly one of $a<0$ or $a>0$ is true, we find the eigenvalues of
the Casimir invariants $C_l$ are given by:

\begin{multline*} 
\chi_{\Lambda}(C_l) = \sum_a (-1)^{[a]} q^{C(\Lambda_0) -
  (\varepsilon_a, \varepsilon_a)} f(a) \Bigl[
  \frac{(q^{(\varepsilon_a, 2 \rho + 2 \Lambda +
  \varepsilon_a)-C(\Lambda_0)}-1)}{(q-q^{-1})} \Bigr]^l \\ \times
  \prod_{b \neq a} \frac{(q^{(\varepsilon_a, 2\rho +2\Lambda +
  \varepsilon_a)} - q^{(\varepsilon_b,
  2\rho+2\Lambda-\varepsilon_b)})} {(q^{(\varepsilon_a, 2\rho
  +2\Lambda+\varepsilon_a)}-q^{(\varepsilon_b, 2\rho
  +2\Lambda+\varepsilon_b)})},
\end{multline*}

\noindent where

\begin{equation*} 
f(a) = \begin{cases} 1 - (q-q^{-1})\frac{q} {(q^{2
(\varepsilon_a,\varepsilon_a + 2\rho + 2\Lambda)} - 1)}, & m=2l, \\ 1
+ (q-q^{-1}) \frac{q^{(\varepsilon_a,\varepsilon_a + 2\rho +
2\Lambda)}} {(q^{2(\varepsilon_a, \varepsilon_a + 2\rho + 2\Lambda)}-
1)}, \quad& a \neq 0, \; m=2l+1, \\ 1, & a=0, \; m=2l+1.
\end{cases}
\end{equation*}

\noindent Throughout we assumed the eigenvalues were distinct.  If
they are not, the calculations are more complicated but the result is
the same.  Thus we have found:

\begin{theorem}

The quantum superalgebra $U_q[osp(m|n)]$, for $m > 2,$ has an infinite family of Casimir invariants
of the form

\begin{equation*}
C_l = (str \otimes I)(\pi(q^{2h_p}) \otimes I) A^l, \qquad l \in
\mathbb{Z}^+,
\end{equation*}

\noindent where

\begin{equation*}
A = \frac{(R^T R - I \otimes I)}{(q-q^{-1})}.
\end{equation*}

\noindent The eigenvalues of the invariants when acting on an
arbitrary irreducible finite-dimensional module with highest weight
$\Lambda$ are given by:

\begin{multline*}
\chi_{\Lambda}(C_l) = \sum_a (-1)^{[a]} q^{C(\Lambda_0) -
  (\varepsilon_a, \varepsilon_a)} f(a) \Bigl[
  \frac{(q^{(\varepsilon_a, \varepsilon_a +2\rho +
  2\Lambda)-C(\Lambda_0)}-1)}{(q-q^{-1})} \Bigr]^l \\ \times \prod_{b
  \neq a} \frac{(q^{(\varepsilon_a, 2\rho +2\Lambda + \varepsilon_a)}
  - q^{(\varepsilon_b, 2\rho+2\Lambda-\varepsilon_b)})}
  {(q^{(\varepsilon_a, 2\rho
  +2\Lambda+\varepsilon_a)}-q^{(\varepsilon_b, 2\rho
  +2\Lambda+\varepsilon_b)})},
\end{multline*}

\noindent where

\begin{equation*}
f(a) = \begin{cases} 1 - (q-q^{-1})\frac{q} {(q^{2
(\varepsilon_a,\varepsilon_a + 2\rho + 2\Lambda)} - 1)}, & m=2l, \\ 1
+ (q-q^{-1}) \frac{q^{(\varepsilon_a,\varepsilon_a + 2\rho +
2\Lambda)}} {(q^{2(\varepsilon_a, \varepsilon_a + 2\rho + 2\Lambda)}-
1)}, \quad& a \neq 0, \; m=2l+1, \\ 1, & a = 0,\; m=2l+1.
\end{cases}
\end{equation*}

\end{theorem}

This completes the calculation of the eigenvalues of an
infinite family of Casimir invariants of $U_q[osp(m|n)]$ when acting
on an arbitrary irreducible highest weight module, provided $m > 2$. 
This had already been done for $U_q[osp(2|n)]$,
using a different method, in \cite{GLZ}. Also every finite-dimensional 
representation of $U_q[osp(1|n)]$ is isomorphic to a
finite-dimensional representation of $U_{-q}[o(n+1)]$ \cite{Zhang},
whose central elements are well-understood.  Hence the eigenvalues of
a family of Casimir invariants, when acting on an arbitrary irreducible
finite-dimensional highest weight module have now been calculated 
for all quantised orthosymplectic superalgebras. Together with the results
for $U_q[gl(m|n)]$ \cite{LinksZhang}, this covers all non-exceptional quantised
superalgebras.


\vspace{0.5cm}
\noindent{\bf Acknowledgements --} We gratefully acknowledge financial
assistance from the Australian Research Council.

\end{document}